\documentclass[a4paper,10pt]{amsart}
\usepackage{amssymb}

\input xypic
\xyoption{all}
\xyoption{arc}

\newtheorem{theo}{Th\'eor\`eme}[section]

\newtheorem{prop}[theo]{Proposition}

\newtheorem{lem}[theo]{Lemme}

\newtheorem{defi}[theo]{D\'efinition}

\def\a{\alpha}
\def\b{\beta}
\def\g{\gamma}

\def\d{\delta}

\def\ph{\varphi}
\def\l{\lambda}

\def\r{\rho}
\def\s{\sigma}

\def\th{\theta}

          \def\rhot{{\tilde{\rho}}}
        \def\sigt{{\tilde{\s}}}

\DeclareMathOperator{\End}{{\mathrm{End}}}

\DeclareMathOperator{\Id}{{\mathrm{Id}}}

\DeclareMathOperator{\INT}{{\mathrm{int}}}

\DeclareMathOperator{\disc}{{\mathrm{disc}}}
\DeclareMathOperator{\Nrd}{{\mathrm{Nrd}}}
\DeclareMathOperator{\Skew}{{\mathrm{Skew}}}
\DeclareMathOperator{\Trd}{{\mathrm{Trd}}}
\DeclareMathOperator{\Sym}{{\mathrm{Sym}}}
\DeclareMathOperator{\Sand}{{\mathrm{Sand}}}

\def\to{\rightarrow}

\def\mapright#1{\hspace{0.2em}\smash{
     \mathop{\rightarrow}\limits^{\SS#1}}\hspace{0.2em}}

\def\vide{\varnothing}

\def\SS{\scriptstyle}

\def\lexp#1#2{\kern\scriptspace\vphantom{#2}^{#1}\kern-\scriptspace#2}

\mathchardef\inferieur="321E
\mathchardef\superieur="321F

\def\eqna{\begin{eqnarray*}}
\def\endeqna{\end{eqnarray*}}

\def\cad{c'est-\`a-dire }

\def\ssi{si et seulement si }

\def\itemth#1{\item[${\mathrm{(#1)}}$]}

\def\Aun{\underline{A}}
\def\Aop{A^{op}}

\begin{document}

\title{Alg\`ebre de Clifford d'un antiautomorphisme}

\author{Anne Cortella}
\address{\noindent 
Labo. de Math. de Besan\c{c}on, UMR CNRS 6623,
Universit\'e de Franche-Comt\'e, 16 Route de Gray, 25030 Besan\c{c}on
Cedex, France} 

\makeatletter
\email{cortella@math.univ-fcomte.fr}

\makeatother


\date{\today}
\def\abstractname{R\'esum\'e}
\begin{abstract} 
Nous d\'efinissons l'alg\`ebre de Clifford d'un antiautomorphisme d'alg\`ebre centrale simple, et la calculons pour les alg\`ebres de degr\'e $2$.

\medskip

\noindent{\sc Abstract.} 
We give a definition of the Clifford algebra of an antiautomorphism of a central simple algebra, and compute it for the algebras of degree $2$.
\end{abstract}

\maketitle

\pagestyle{myheadings}

\markboth{\sc A. Cortella}{\sc Alg\`ebre de Clifford}

L'alg\`ebre de Clifford d'une involution lin\'eaire de type orthogonal sur une alg\`ebre centrale simple $A$ a \'et\'e d\'efinie par Jacobson dans \cite{J}, par descente galoisienne, pour g\'en\'eraliser la partie paire de l'alg\`ebre de Clifford d'une forme quadratique. Tits \cite{TI} l'a ensuite d\'efinie comme un quotient de l'alg\`ebre tensorielle de l'espace vectoriel sous-jacent \`a $A$.

Dans cet article, une d\'efinition de l'alg\`ebre de Clifford pour un antiautomorphisme $\s$, lin\'eaire, non n\'ecessairement involutif, sur une alg\`ebre centrale simple $A$, est propos\'ee (d\'efinition \ref{clifford}). Elle g\'en\'eralise la d\'efinition de Tits. Cela d\'efinit un nouvel invariant de la classe d'isomorphie de $(A,\s)$ (th\'eor\`eme \ref{invt}), se comportant bien par extension des scalaires. Pour $A$ d\'eploy\'ee, c'est la partie paire de ce qu'on peut d\'efinir comme la g\'en\'eralisation aux formes bilin\'eaires non necessairement sym\'etriques de l'alg\`ebre de clifford classique des formes quadratiques (th\'eor\`eme \ref{deploye} et d\'efinition \ref{clif bilin}).

Avant de d\'etailler cette d\'efinition et ces propri\'et\'es dans la deuxi\`eme partie, nous commencerons, dans la premi\`ere partie, par rappeler les d\'efinitions du discriminant et de l'alg\`ebre de Clifford dans le cas involutif orthogonal, puis celles d'invariants des antiautomorphismes issus de \cite{CT} et utilis\'es dans notre d\'efinition principale.

Enfin, une troisi\`eme partie est consacr\'ee au calcul de notre nouvel invariant si l'alg\`ebre $A$ est de degr\'e $2$. On fait alors le lien avec le discriminant et l'asym\'etrie de l'antiautomorphisme (proposition \ref{deg 2}).

\bigskip

\noindent{\sc Notations - } Soit $F$ 
un corps de caract\'eristique diff\'erente de $2$. Si $V$ et $W$ sont deux $F$-espaces vectoriels, si 
$\psi : V \to W$ est une application $F$-lin\'eaire et si $L$ est une extension 
de $F$, on pose $V_L = V \otimes_F L$ et on note $\psi_L : V_L \to W_L$ 
l'application $L$-lin\'eaire $\psi \otimes_F \Id_L$. Si $\s$ est un 
$F$-endomorphisme de $V$, on pose $\Sym(V,\s)=\{v \in V~|~\s(v)=v\}$ 
et $\Skew(V,\s)=\{v \in V~|~\s(v)=-v\}$. L'alg\`ebre tensorielle de 
$V$ sur $F$ sera not\'ee $T(V)$.

Si $A$ est une $F$-alg\`ebre centrale simple, on note $\Aun$ l'espace vectoriel sous-jacent \`a $A$, $A^{op}$ l'alg\`ebre oppos\'ee. On note respectivement $\Trd : A \to F$ et $\Nrd : A \to F$ la trace r\'eduite et la norme r\'eduite. 


\bigskip

\section{Quelques rappels}

\subsection{Le discriminant et l'alg\`ebre de Clifford d'une involution 
de type orthogonal} 
Ici, $A$ est une alg\`ebre centrale simple sur $F$ et $\s$ est une involution de type orthogonal
sur $A$, \cad un antiautomorphisme $F$-lin\'eaire involutif de l'alg\`ebre $A$ tel que~: si $L$ est un corps qui d\'eploie $A$, en identifiant $A_L \simeq \End_L(L^n)$, alors 
$\s_L$ est l'adjonction $\s_b$ pour un 
forme bilin\'eaire sym\'etrique non d\'eg\'en\'er\'ee $b$ sur $V=L^n$~:
\begin{equation}\label{adjonction}
\forall~f \in \End_L V,~\forall~x,y \in V,~b(f(x),y)=b(x,\s_b(f)(y)).
\end{equation}

Si $A$ est d\'eploy\'ee, la classe de similitude de $b$ est un invariant de la classe d'isomorphie de $(A,\s)$.
\begin{prop}[Knus-Parimala-Sridharan {\cite[proposition 7.1, page 81]{BOI}}]\label{PKS}
Si $\deg A=2m$ est pair, alors $\Skew(A,\s) \cap A^\times \not= \vide$ et 
$\Nrd x \in F^\times/(F^\times)^2$ est ind\'ependant de $x$ 
pris dans cet ensemble.
\end{prop}

Cela leur permet de d\'efinir le discriminant 
pour les alg\`ebres de degr\'e pair $2m$~:

\begin{defi}\label{PKS defi}
On appelle discriminant de $(A,\s)$, et on note $\disc \s$, 
la classe de $(-1)^m \Nrd x$ dans $F^\times/(F^\times)^2$, o\`u $x$ est un \'el\'ement 
quelconque de $\Skew(A,\s) \cap A^\times$.
\end{defi}

Cela d\'efinit un invariant de la classe d'isomorphie de $(A,\s)$ qui \'etend 
la notion de discriminant (\`a signe) des formes bilin\'eaires sym\'etriques~: 
si $b$ est une telle forme sur un espace vectoriel $V$ de dimension $n$, 
on d\'efinit $\disc b=(-1)^{\frac{n(n-1)}{2}} \det b \in F^\times/(F^\times)^2$, qui est un invariant de la classe de similitude de $b$ ; si de plus $n=2m$, alors $\disc b=\disc \s_b$. 
De plus, le discriminant se comporte bien pour l'extension des scalaires 
\cite[proposition 7.3, page 81]{BOI}.

\medskip

L'alg\`ebre de Clifford d'une involution de type orthogonal 
a \'et\'e \`a l'origine d\'efinie par Jacobson \cite{J}, et c'est le centre de cette alg\`ebre qui a donn\'e la premi\`ere d\'efinition du discriminant.
Tits \cite{TI} en a cependant donn\'e une d\'efinition plus facile \`a mettre en \oe uvre. 
Pour cela, il utilise le sandwich.

L'application $\Sand : A \otimes A^{op} \to \End_F \Aun$, 
$x \otimes y \mapsto (a \mapsto xay)$ est un isomorphisme d'alg\`ebres 
appel\'e sandwich. Si $\s$ est une involution orthogonale sur $A$ et si $u \in A \otimes A^{op}$, 
alors $a \mapsto (\Sand u)(\s(a))$ est un \'el\'ement de $\End_F \Aun$, 
donc il existe $\s_2(u) \in A \otimes A^{op}$ tel que, pour tout 
$a \in A$, $(\Sand \s_2(u))(a)=(\Sand u)(\s(a))$. Cela d\'efinit 
une application lin\'eaire involutive $\s_2$ sur l'espace vectoriel $\Aun \otimes \Aun$. 
Notons encore $\mu$ la multiplication : $\Aun \otimes \Aun \to A$ ; $u \mapsto (\Sand u)(1)$.

\begin{defi}[Tits {\cite[definition 8.7, page 92]{BOI}}]\label{Tits defi}
L'alg\`ebre de Clifford $C(A,\s)$ est le quotient 
$$C(A,\s)= T(\Aun)/(J_1(A,\s) + J_2(A,\s))$$
o\`u 
\begin{itemize}
\itemth{1} $J_1(A,\s)$ est l'id\'eal de $T(\Aun)$ engendr\'e par les 
$s-\frac{1}{2}\Trd s$ pour $s \in \Sym(\Aun,\s)$. 

\itemth{2} $J_2(A,\s)$ est l'id\'eal de $T(\Aun)$ engendr\'e par 
les $u-\frac{1}{2} \mu(u)$ pour $u \in \Sym(\Aun \otimes \Aun,\s_2)$.
\end{itemize}
\end{defi}

Cela d\'efinit encore un invariant de la classe d'isomorphie de $(A,\s)$. 
En particulier, si $A$ est d\'eploy\'ee et $\s=\s_b$ avec $b$ bilin\'eaire 
sym\'etrique sur $V$, alors, en identifiant via 
$\ph_b : V \otimes V \mapright{\sim} \Aun=\underline{\End_F V}$ ; 
$v \otimes w \mapsto (x \mapsto vb(w,x))$, on trouve 
$C(A,\s)=T(V \otimes V)/(I_1 + I_2)$, o\`u $I_1$ est l'id\'eal engendr\'e par 
les $v \otimes v - \frac{1}{2} b(v,v)$, $v \in V$, 
et $I_2$ par les $u \otimes v \otimes v \otimes w - \frac{1}{2} b(v,v) (u \otimes w)$, 
$u$, $v$, $w \in V$, et donc $C(A,\s)$ est la partie paire de l'alg\`ebre gradu\'ee 
$T(V)/<v \otimes v - b(v,v), v \in V>$, qui n'est autre que l'alg\`ebre de 
Clifford classique pour une forme quadratique (ceci ne d\'ependant 
que de la classe de similitude de $b$, on peut enlever le facteur $\frac{1}{2}$). 

De plus, $C(A,\s)$ se comporte bien par extension des scalaires. 

\smallskip

Signalons encore une propri\'et\'e structurelle importante.

\begin{theo}[{\cite[th\'eor\`eme 8.10, page 94]{BOI}}]\label{structure}
Si $\deg A = 2m$, le centre de $C(A,\s)$ est l'alg\`ebre \'etale 
$Z=F[X]/(X^2-\disc \s)$. Si c'est un corps, alors $C(A,\s)$ est une $Z$-alg\`ebre 
centrale simple de degr\'e $2^{m-1}$. Sinon, $Z \simeq F \times F$ 
et $C(A,\s)$ est le produit de deux $F$-alg\`ebres centrales simples de 
degr\'e $2^{m-1}$.
\end{theo}

Enfin, $\s$ induit sur $T(\Aun)$ une involution par~: 
$\underline{\s}(a_1 \otimes \dots \otimes a_r) = \s(a_r) \otimes \dots \otimes \s(a_1)$, 
qui passe au quotient et d\'efinit ainsi une involution dite canonique 
sur $C(A,\s)$.

\subsection{L'asym\'etrie et le discriminant d'un antiautomorphisme} 
Soit maintenant une alg\`ebre centrale simple $A$ sur le corps $F$ munie 
d'un antiautomorphisme $F$-lin\'eaire $\s$ non n\'ecessairement involutif. Par extension des 
scalaires \`a un corps $L$ scindant $A$ et en identifiant $A_L \simeq \End_L V$, 
o\`u $V$ est un $L$-espace vectoriel, l'antiautomorphisme $\s_L$ est l'adjonction pour une forme bilin\'eaire non d\'eg\'en\'er\'ee $b : V \times V \to L$ non 
n\'ecessairement sym\'etrique, d\'efinie \`a similitude pr\`es par la classe d'isomorphie de $\s$, \cad que $\s_L$ v\'erifie 
\ref{adjonction}.

\begin{prop}[{\cite[propositions 3 et 4]{CT}}]

\hspace{1cm}

\begin{itemize}

\medskip

\itemth{1} Si $b$ est une forme bilin\'eaire non d\'eg\'en\'er\'ee 
sur $V$, il existe une unique application lin\'eaire $\g_b$ sur 
$\underline{\End_F V}$ satisfaisant 
$$\forall~x,y \in V,~\forall~f \in \End_F V,~b(x,f(y))=b(y,\g_b(f)(x)).$$

\itemth{2} Si $\s$ est un antiautomorphisme sur $A$ ($F$-lin\'eaire), il existe 
une unique application lin\'eaire $\g_\s$ sur $\Aun$ telle que 
si $L$ scinde $A$, si $\th$ est un isomorphisme de 
$A_L=A \otimes_F L \mapright{\sim} \End_L V$ et si $b$ est une forme bilin\'eaire 
sur $V$ telle que $\th \circ \s_L \circ \th^{-1} = \s_b$, alors 
$\th \circ (\g_\s \otimes_F \Id_L) \circ \th^{-1} = \g_b$.
\end{itemize}
De plus, pour tous $x$, $y$, $z \in A$, $\g_\s(xyz)=\s(z)\g_\s(y)\s^{-1}(x)$  
et $\g_\s^2 = \Id_A$. 
\end{prop}

\begin{defi}
L'asym\'etrie de $(A,\s)$ est l'\'el\'ement $a_\s=\g_\s(1) \in A^\times$. 
\end{defi}

Cet \'el\'ement v\'erifie en particulier~:
\begin{itemize}
\itemth{1} Si $\s=\s_b$, alors, pour tous $x$, $y \in V$, $b(x,y)=b(y,a_\s(x))$. 

\itemth{2} Pour tout $\a \in A$, $\g_\s(\a)=\s(\a) a_\s$.

\itemth{3} $\s^2=\INT a_\s$. 

\itemth{4} $\s(a_\s)=a_\s^{-1}$.

\itemth{5} $\s$ est une involution orthogonale \ssi $a_\s=1$, et 
$\s$ est une involution symplectique \ssi $a_\s=-1$.
\end{itemize}
De plus, la classe d'isomorphie du couple $(\Aun,\g_\s)$ et la classe 
de conjugaison de $a_\s$ sont des invariants de la classe d'isomorphie de 
$(A,\s)$. 

\smallskip

Nous pouvons alors, si $\deg A=2m$, \'etendre la d\'efinition du 
discriminant de $(A,\s)$ aux antiautomorphismes en utilisant $\g_\s$. 

\begin{prop}[{\cite[lemme 3]{CT}}]\label{CT}
Si $\deg A=2m$ est pair, alors $\Skew(\Aun,\g_\s) \cap A^\times \not= \vide$ et 
$\Nrd x \in F^\times/(F^\times)^2$ est ind\'ependant de $x$ 
pris dans cet ensemble.
\end{prop}

\begin{defi}
On appelle discriminant de $(A,\s)$, et on note $\disc \s$, 
la classe de $(-1)^m \Nrd x$ dans $F^\times/(F^\times)^2$, o\`u $x$ est un \'el\'ement 
quelconque de $\Skew(\Aun,\g_\s) \cap A^\times$.
\end{defi}

En particulier, si $\s=\s_b$, $\disc \s = (-1)^{\frac{n(n-1)}{2}} \det b$ et, 
si $\s$ est une involution orthogonale, le discriminant est bien celui 
de la definition \ref{PKS defi}.
Par ailleurs, si $1-a_\s$ est inversible, alors $\disc \s=(-1)^m\Nrd(1-a_\s)$.

\section{L'alg\`ebre de Clifford d'un antiautomorphisme}

\subsection{D\'efinition et invariance}

Fixons une $F$-alg\`ebre centrale simple $A$, munie d'un antiautomorphisme 
$F$-lin\'eaire $\s$. 
Nous d\'efinissons ici l'alg\`ebre de Clifford de $(A,\s)$ 
en suivant la d\'efinition \ref{Tits defi} de Tits. 

\begin{defi}\label{clifford}
L'alg\`ebre de Clifford $C(A,\s)$ est le quotient 
$$C(A,\s)= T(\underline{A})/(J_1(A,\s) + J_2(A,\s))$$
o\`u 
\begin{itemize}
\itemth{1} $J_1(A,\s)$ est l'id\'eal de $T(\Aun)$ engendr\'e par les 
$s-\frac{1}{2}\Trd s$ pour $s \in \Sym(\Aun,\g_\s).$ 

\itemth{2} $J_2(A,\s)$ est l'id\'eal de $T(\Aun)$ engendr\'e par 
les $u-\frac{1}{2} \mu_\s(u)$ pour $u \in \Sym(\Aun \otimes \Aun,\g_{\sigt,2})$, avec 

$\bullet$ $\sigt$ est l'antiautomorphisme $\sigt = (\INT a_\s) \circ \s$ de $A$ et donc $\g_{\sigt}$ est l'application lin\'eaire involutive de $\Aun$ d\'efinie par 
$\g_{\sigt}(x)=a_\s \g_\s(x) a_\s$ ; alors $\g_{\sigt,2}$ est l'application lin\'eaire involutive induite par $\g_{\sigt}$ sur $\Aun \otimes \Aun$ gr\^ace au sandwich par 
$$\forall u\in A\otimes\Aop\quad\forall x\in A\quad (\Sand \g_{\sigt,2}(u))(x)=(\Sand u)(\g_{\sigt}(x))\,;$$

$\bullet$ si $u\in\Aun\otimes\Aun,\quad\mu_\s(u)=(\Sand u)(a_\s)$.
\end{itemize}
\end{defi}

Avec cette d\'efinition, il est clair que si $\s$ est une involution 
de type orthogonal, alors $C(A,\s)$ n'est autre que l'alg\`ebre 
de Clifford de $(A,\s)$ d\'efinie par Tits. 

\begin{theo}\label{invt}
La classe d'isomorphie de $C(A,\s)$ est un invariant de la classe 
d'isomorphie de $(A,\s)$.
\end{theo}

\begin{proof}
Supposons $(A,\s)$ et $(A',\r)$ isomorphes, \cad qu'il existe un isomorphisme 
$\d : A \mapright{\sim} A'$ tel que $\r= \d \circ \s \circ \d^{-1}$. 
On peut alors supposer que $A=A'$ et, comme elle est centrale simple, 
d'apr\`es le th\'eor\`eme de Skolem-Noether, $\d$ est int\'erieur. 
Il existe donc $w \in A^\times$ tel que $\d = \INT w$. Ainsi, 
$\r=\INT(\r(w)w) \circ \s$. Nous allons alors relier $a_\s$ \`a $a_\r$, 
$\g_\s$ \`a $\g_\r$ et $\g_{\s,2}$ \`a $\g_{\r,2}$. 

\begin{prop}\label{conjug}
Si $\s$ et $\r$ sont deux antiautomorphismes isomorphes de $A$ et si 
$w \in A^\times$ v\'erifie $\r \circ \INT w = \INT w \circ \s$, alors
\begin{itemize}
\itemth{1} $a_\r=w a_\s w^{-1}$.

\itemth{2} $\g_\r \circ \INT w = \INT w \circ \g_\s$.

\itemth{3} $\g_{\r,2} \circ \INT(w \otimes w^{-1}) = 
\INT(w \otimes w^{-1}) \circ \g_{\s,2}$. 
\end{itemize}
\end{prop}

\begin{proof}
(1) est un calcul direct \`a partir de \cite[proposition 7]{CT}~: si 
$\r=(\INT u) \circ \s$, alors $a_\r = u \s(u)^{-1} a_\s$ que l'on 
applique \`a $u=\r(w)w=w\s(w)$. 

(2) est alors obtenu gr\^ace \`a $\g_\s(\a)=\s(\a) a_\s$ si $\a \in A$. 

Pour (3), on utilise la structure d'alg\`ebre donn\'ee sur $\Aun \otimes \Aun$ 
par $A \otimes A^{op}$, qui fait du sandwich un isomorphisme d'alg\`ebres. On peut alors d\'emontrer :

\begin{lem}\label{lemme conjug}
Si $w\in A$ et $u\in A\otimes\Aop$, alors $\Sand(\INT(w\otimes w^{-1})(u))=\INT w\circ\Sand u\circ\INT w^{-1}$.

\end{lem}

\noindent En effet, si $u=\a \otimes \b \in A \otimes A^{op}$ et $x\in A$, alors

\eqna
[\Sand(\INT(w\otimes w^{-1})(\a \otimes \b))](x)&=& [\Sand((w\otimes w^{-1})(\a \otimes \b)(w\otimes w^{-1})^{-1})](x)\\
&=&[\Sand(w\a w^{-1}\otimes w\b w^{-1})](x)\\
&=&w\a w^{-1}xw\b w^{-1}\\
&=&w[(\Sand(\a\otimes\b))(w^{-1}xw)]w^{-1},
\endeqna                                   
d'o\`u le lemme.

\medskip

De ce lemme on d\'eduit que si $w$ est comme dans l'\'enonc\'e de la proposition et $u\in A\otimes\Aop$, alors

\eqna
&&\Sand[\INT(w\otimes w^{-1})\circ\g_{\s,2}\circ\INT(w\otimes w^{-1})^{-1}(u)]\\
&=&\INT w\circ\Sand[\g_{\s,2}(\INT(w^{-1}\otimes w)(u))]\circ\INT w^{-1}\\
&=&\INT w\circ\Sand(\INT(w^{-1}\otimes w)(u))\circ\g_\s\circ\INT w^{-1}\\
&=&\INT w\circ\INT w^{-1}\circ\Sand(u)\circ\INT w\circ\g_\s\circ\INT w^{-1}\\
&=&\Sand u\circ\g_\r=\Sand(\g_{\r,2}(u)).
\endeqna

Donc finalement 
$$\g_{\r,2}=\INT(w\otimes w^{-1})\circ\g_{\s,2} \circ \INT(w\otimes w^{-1})^{-1},$$  comme annonc\'e.          \end{proof}

D'apr\`es cette proposition, $\Sym(A,\g_\r)=(\INT w)(\Sym(A,\g_\s))$ et donc, 
en notant $T(w)$ l'automorphisme de l'alg\`ebre $T(\Aun)$ 
induit par l'application lin\'eaire $\INT w$ sur $\Aun$, 
$J_1(A,\r)=T(w)(J_1(A,\s))$ : en effet, si $s \in \Sym(A,\g_\s)$, 
alors $T(w)(s -\frac{1}{2} \Trd s)=wsw^{-1} - \frac{1}{2} \Trd(wsw^{-1})$. 

\medskip

Remarquons maintenant que $\sigt=(\INT a_\s) \circ \s$ et 
$\rhot=(\INT a_\r) \circ \r= (\INT w) \circ \sigt \circ (\INT w)^{-1}$. 
Donc d'apr\`es le (3) du lemme pr\'ec\'edent, 
$\g_{\rhot,2} \circ (\INT w \otimes w^{-1}) = 
(\INT w \otimes w^{-1}) \circ \g_{\sigt,2}$, 
ce qui redonne pour les \'el\'ements sym\'etriques 
$\Sym(\Aun \otimes \Aun,\g_{\rhot,2})=
(\INT w \otimes w^{-1})(\Sym(\Aun \otimes \Aun,\g_{\sigt,2}))$. 
Il ne reste plus qu'\`a remarquer que $\INT w \otimes w^{-1}$ est la 
restriction de $T(w)$ \`a $\Aun \otimes \Aun$ et que, sur $\Aun \otimes \Aun$, 
on a l'\'egalit\'e 
$$\mu_\r \circ T(w) = T(w) \circ \mu_\s$$ d'apr\`es le lemme \ref{lemme conjug}.

On en d\'eduit que $J_2(A,\r)=T(w)(J_2(A,\s))$, ce qui ach\`eve la d\'emonstration
du th\'eor\`eme.                                    
  
\end{proof}

\subsection{Le cas d\'eploy\'e}

\begin{theo}\label{deploye}
L'alg\`ebre de Clifford est compatible \`a l'extension des 
scalaires : si $L$ est une extension de $F$, alors $C(A_L,\s_L)=C(A,\s)\otimes L$. De plus, si $A$ est d\'eploy\'ee et si $\s=\s_b$, o\`u 
$b$ est une forme bilin\'eaire sur l'espace vectoriel $V$ sur $F$, 
alors $C(A,\s)$ est la partie paire $C_0(V,b)$ de l'alg\`ebre 
$C(V,b)=T(V)/<a_b(v) \otimes v - b(v,v),v \in V>$. 
\end{theo}

C'est cette deuxi\`eme partie du th\'eor\`eme qui justifie 
que l'on prenne dans la d\'efinition de $J_2$ les \'el\'ements 
sym\'etriques pour $\g_{\sigt,2}$ au lieu de prendre tout simplement 
ceux sym\'etriques pour $\g_{\s,2}$ : m\^eme si cette autre d\'efinition 
donnerait aussi, en utilisant la proposition \ref{conjug}, un invariant de la classe d'isomorphie de $(A,\s)$, il ne serait pas en g\'en\'eral dans le cas d\'eploy\'e la partie paire d'un quotient de $T(V)$. 

Ce th\'eor\`eme conduit \`a la d\'efinition suivante :

\begin{defi}\label{clif bilin}
Soit $b$ une forme bilin\'eaire non d\'eg\'en\'er\'ee sur l'espace vectoriel $V$ sur $F$. On appelle alg\`ebre de clifford de $b$ l'alg\`ebre $C(V,b)=T(V)/<a_b(v) \otimes v - b(v,v),v \in V>$.
\end{defi}

Cette alg\`ebre n'est en g\'en\'eral pas du type de celles d\'efinie par Hannabus dans \cite{H}, en particulier elles diff\`erent d\`es que $b$ est antisym\'etrique, ou plus g\'en\'eralement d\`es que l'asym\'etrie de $b$ admet $-1$ comme valeur propre.

\begin{proof}
La premi\`ere assertion est claire. Montrons la deuxi\`eme~: 
supposons que $A$ est d\'eploy\'ee, $A=\End_F V$ et 
$\s=\s_b$ o\`u $b : V \times V \to F$ est bilin\'eaire non d\'eg\'en\'er\'ee. 
Pour simplifier, on oubliera ici d'indexer $a$, $\g$ et $\ph$ par $\s$. 
Rappelons que $\ph$ est l'isomorphisme $V \otimes V \mapright{\sim} \Aun$ ; 
$v \otimes w \mapsto (x \mapsto vb(w,x))$. Il induit un 
isomorphisme d'alg\`ebres $T(\ph) : T(V \otimes V) \mapright{\sim} T(\Aun)$. 
Nous devons alors calculer $T(\ph)^{-1}(J_1)$ et $T(\ph)^{-1}(J_2)$. 

\medskip

Pour calculer $T(\ph)^{-1}(J_1)$, nous allons utiliser le lemme calculatoire suivant :
 
\begin{lem}\label{lemmeJ_1}
Soient $v,w\in V$. Alors 
$$\s\circ\ph(v\otimes w)=\ph(a(w)\otimes v),\quad \ph(v\otimes w)=\ph(v\otimes a(w))\circ a\quad\mbox{et}\quad a\circ\ph(v\otimes w)=\ph(a(v)\otimes w),$$ dont on d\'eduit que $\g\circ\ph(a(v)\otimes w)=\ph(a(w)\otimes v)$.

\end{lem}

\noindent En effet, comme $\s =\s_b$, il v\'erifie 
$\forall~f \in \End_L(V),~\forall~x,y \in V,~b(f(x),y)=b(x,\s(f)(y)).$
En particulier, si $f=\ph(v \otimes w)$, alors pour $x,y\in V$,
\eqna
b(x,\s(\ph(v \otimes w))(y))&=& b(\ph(v \otimes w)(x),y) \\
&=& b(v b(w,x),y) \\
&=& b(v,y) b(w,x) \\
&=& b(x,a(w)) b(v,y) \\
&=& b(x,\ph(a(w) \otimes v)(y))
\endeqna
d'o\`u la premi\`ere \'egalit\'e. De m\^eme 
\eqna
b(\ph(v \otimes a(w))\circ a(x),y)&=& b(vb(a(w),a(x)),y) \\
&=& b(v,y)b(a(w),a(x)) \\
&=& b(v,y) b(w,x) \\
&=& b(\ph(v\otimes w)(x),y)
\endeqna
ce qui prouve la deuxi\`eme \'egalit\'e. La troisi\`eme se montre de la m\^eme mani\`ere. On d\'eduit de ces trois \'egalit\'es que $\g\circ\ph(a(v)\otimes w)=\s\circ\ph(a(v)\otimes w)\circ a=\ph(a(w)\otimes a(v))\circ a=\ph(a(w)\otimes v),$ ce qui ach\`eve la preuve du lemme.

\medskip

Soit maintenant $f\in\End V$, elle est sym\'etrique pour $\g$ si et seulement si $f=\displaystyle\frac{1}{2}\bigl(f+\g(f)\bigr)$. \'Ecrivons $f$ comme l'image par $\ph$ d'une combinaison lin\'eaire de tenseurs \'el\'ementaires : 

\noindent $f=\ph(\sum\limits_{i=1}^n\l_{i}a(v_i)\otimes w_i)$. Alors d'apr\`es le lemme $f+\g(f)=\sum\limits_{i=1}^n\l_{i}(\ph(a(v_i)\otimes w_i)+\ph(a(w_i)\otimes v_i))$.

On en d\'eduit que $f \in \Sym(A,\g)$ 
si et seulement si elle est combinaison lin\'eaire de termes de la forme $\ph(a(v)\otimes w)+\ph(a(w)\otimes v)$. Or ceci n'est autre que $\ph(a(v+w)\otimes (v+w))-\ph(a(v)\otimes v)-\ph(a(w)\otimes w)$. Ainsi $f$ est sym\'etrique si et seulement si $\ph^{-1}(f)$ est combinaison lin\'eaire
d'\'el\'ements 
du type $a(w) \otimes w$ avec $w \in V$. 
De plus, si $v$, $w \in V$, $\ph(v \otimes w)$ a pour trace $b(w,v)=b(v,a(w))$, donc 
$\Trd(\ph(a(w) \otimes w))=b(a(w),a(w))=b(w,w)$. Ainsi $T(\ph)^{-1}(J_1)$ 
est l'id\'eal $I=<a(w) \otimes w - \frac{1}{2} b(w,w), w \in V>$. 

\medskip

Pour d\'eterminer  $T(\ph)^{-1}(J_2)$, nous utiliserons encore un lemme calculatoire :
 
\begin{lem}\label{lemmeJ_2}
Soient $y,z,v,w,s,t\in V$. Alors 
\eqna
\bigl(\Sand(\ph(v\otimes w)\otimes\ph(s\otimes t))\bigr)(\ph(z\otimes y))&=&\bigl(\Sand(\ph(v\otimes a^{-1}(s))\otimes\ph(a(w)\otimes t))\bigr)(\ph(y\otimes z))\quad\mbox{et}\\
\g_{\sigt,2}(\ph(v\otimes w)\otimes\ph(s\otimes t))&=&\ph(v\otimes a(s))\otimes\ph(a^{-1}(w)\otimes t).
\endeqna

\end{lem}

\smallskip

\noindent En effet, si de plus $x\in V$, alors
\eqna
\bigl(\Sand(\ph(v\otimes w)\otimes\ph(s\otimes t))\bigr)(\ph(z\otimes y))(x)
&=&\ph(v\otimes w)\circ\ph(z\otimes y)\circ\ph(s\otimes t)(x)\\
&=&vb(w,z)b(y,s)b(t,x)\\
&=&vb(a^{-1}(s),y)b(z,a(w))b(t,x)\\
&=&\ph(v\otimes a^{-1}(s)\circ\ph(y\otimes z)\circ\ph(a(w)\otimes t)(x)
\endeqna
d'o\`u la premi\`ere \'egalit\'e. Alors si $u \in A \otimes A^{op}$, 
$\g_{\sigt,2}(u)$ est d\'efini par~: si $y$, $z \in V$, 
\eqna
\Sand(\g_{\sigt,2}(u))(\ph(y \otimes z)) & =& 
(\Sand u)(\g_{\sigt}(\ph(y \otimes z)) \\
&=& (\Sand u)\bigl(a.(\s\circ\ph)(y \otimes z). a^2\bigr) \\
&=& (\Sand u)\bigl((a.\ph(a(z) \otimes y).a^2)\bigr)\\
&=& (\Sand u)\bigl((a^2.\ph(z\otimes y).a^2)\bigr)\\
&=& (\Sand (u.a^2\otimes a^2))(\ph(z\otimes y))
\endeqna
d'apr\`es le lemme \ref{lemmeJ_1}.
Si de plus $u=\ph(v \otimes w)\otimes\ph(s \otimes t)$, toujours d'apr\`es ce lemme, 

\eqna
(\Sand\g_{\sigt,2}(u))\ph(y,z)&=&\bigl(\Sand(\ph(v\otimes w)a^2\otimes a^2\ph(s\otimes t))\bigr)(\ph(z\otimes y))\\
&=&\bigl(\Sand(\ph(v\otimes a^{-2}(w))\otimes\ph(a^2(s)\otimes t))\bigr)(\ph(z\otimes y))\\
&=&\bigl(\Sand(\ph(v\otimes a(s))\otimes\ph(a^{-1}(w)\otimes t))\bigr)(\ph(y\otimes z))
\endeqna
en utilisant la premi\`ere \'egalit\'e. D'o\`u le r\'esultat.

\medskip

Un \'el\'ement $u=\sum\limits_{i,j,k,l=1}^n\l_{ijkl}\ph(e_i\otimes a(e_j))\otimes\ph(e_k\otimes a(e_l))\in\Aun\otimes\Aun$ est donc sym\'etrique pour $\g_{\sigt,2}$ si et seulement si 
$\forall i,l\in\{1,...,n\}\quad \sum\limits_{j,k=1}^n\l_{ijkl}a(e_k)\otimes e_j =\sum\limits_{j,k=1}^n\l_{ijkl}a(e_j)\otimes e_k$, c'est-\`a-dire, en appliquant $1\otimes a$, si et seulement si les vecteurs $\sum\limits_{j,k=1}^n\l_{ijkl}a(e_j)\otimes e_k$, sont tous sym\'etriques pour $\g$.

On en d\'eduit que $u \in \Sym(\Aun \otimes \Aun, \g_{\sigt,2})$ 
si et seulement si $T(\ph)^{-1}(u)$ est combinaison lin\'eaire 
de tenseurs de la forme $(v\otimes a(s))\otimes(s\otimes t)$ avec 
$v$, $s$ et $t \in V$. Or, pour $u=\ph(v\otimes a(s))\otimes \ph(s\otimes t)$ 
et $x \in V$, 
\eqna
\mu_\s(u)(x)&=& \bigl((\Sand u)(a)\bigr)(x)\\
&=& \bigl(\ph(v\otimes a(s))\circ a\circ\ph(s\otimes t)\bigr)(x) \\
&=& \bigl(\ph(v\otimes a(s))\circ\ph(a(s)\otimes t)\bigr)(x) \\
&=& v . b(a(s),a(s)). b(t,x) \\
&=& v. b(s,s).b(t,x) \\
&=& \ph(b(s,s)(v \otimes t))(x).
\endeqna
Donc $T(\ph)^{-1}(J_2)=
<v \otimes a(s) \otimes s\otimes t -\frac{1}{2} b(s,s) v \otimes t, v,s,t \in V>$. 
Ceci ach\`eve de d\'emontrer le r\'esultat, quitte \`a prendre la forme 
similaire $\frac{1}{2}b$ \`a la place de $b$.
\end{proof}

\subsection{Dimension finie de l'alg\`ebre de Clifford}

\begin{theo}\label{dimension}

L'alg\`ebre de Clifford d'une forme bilin\'eaire non d\'eg\'en\'er\'ee sur un espace vectoriel $V$ est de dimension finie inf\'erieure ou \'egale \`a 
$2^{\dim V}$.

L'alg\`ebre de Clifford d'un antiautomorphisme d'une alg\`ebre centrale simple $A$ est de dimension finie inf\'erieure ou \'egale \`a $2^{\deg A-1}$.
\end{theo}

\begin{proof}
D'apr\`es le th\'eor\`eme \ref{deploye}, il suffit de montrer la premi\`ere partie de l'\'enonc\'e pour une forme bilin\'eaire non d\'eg\'en\'er\'ee $b$ d'asym\'etrie trigonalisable sur le corps $F$. Pla\c cons nous alors dans une base $(e_i)_{i=1,...,n}$ de $V$ dans laquelle la matrice de $a=a_b$ est la matrice triangulaire sup\'erieure $(\a_{ij})_{i,j=1,...n}$. Notons $(b_{ij})$ la matrice de la forme $b$ dans cette base.

Alors $C(V,b)=T(V)/I$ o\`u $I$ est l'id\'eal engendr\'e par les \'el\'ements $a(e_i)\otimes e_i-b_{ii}$ pour $i\in\{1,...n\}$ et $a(e_i)\otimes e_j+a(e_j)\otimes e_i-b_{ij}-b_{ji}$ pour $i,j\in\{1,...,n\}$ tels que $i<j$.

En notant $u\times v$ le produit dans $C(V,b)$ des classes des \'el\'ements $u,v\in V$, on obtient ainsi si $i,j$ sont comme ci-dessus :

$$\displaystyle\sum\limits_{k=1}^i\a_{ki}e_k\times e_i=b_{ii}\quad\mbox{et}\quad \displaystyle\sum\limits_{k=1}^i\a_{ki}e_k\times e_j +\sum\limits_{l=1}^j\a_{lj}e_l\times e_i=b_{ij}+b_{ji}.
$$

Or $a$ est inversible donc les $\a_{ii}$ sont tous non nuls, ce qui permet d'\'ecrire

$$\displaystyle{e_i\times e_i =\frac{b_{ii}}{\a_{ii}}-\sum\limits_{k=1}^{i-1}\frac{\a_{ki}}{\a_{ii}}e_k\times e_i           \quad\mbox{et}\quad          e_j\times e_i =\frac{b_{ij}+b_{ji}}{\a_{jj}}-\sum\limits_{l=1}^{j-1}\frac{\a_{lj}}{\a_{jj}}e_l\times e_i -\sum\limits_{k=1}^{i}\frac{\a_{ki}}{\a_{jj}}e_k\times e_j}.$$

Cette derni\`ere \'egalit\'e signifie que si $i<j$, alors $e_j\times e_i$ est combinaison lin\'eaire de $1$, des $e_l\times e_i$ pour $l<j$ et des $e_k\times e_j$ pour $k\leq i<j$. En appliquant alors successivement ce r\'esultat \`a chacun des $e_l\times e_i$ pour $i<l<j$, on prouve que $e_j\times e_i$ est combinaison lin\'eaire de $1$ et des $e_k\times e_l$ pour $k\leq i\leq l\leq j$ : nous appellerons cela des formules de commutation.

Revenons alors \`a la premi\`ere \'egalit\'e. Elle implique que $e_ i\times e_i$ est combinaison lin\'eaire de $1$ et des $e_k\times e_i$ pour $k<i$. Donc finalement, si $i_1,...,i_k\in\{1,...n\}$, le produit $e_{i_1}\times...\times e_{i_k}$ dans $C(V,b)$ peut s'\'ecrire comme combinaison lin\'eaire de $1$ et des produits 
$e_{j_1}\times...\times e_{j_l}$ pour tous les $l$-uplets $j_1,...j_l\in\{1,...n\}$ satisfaisant $j_1<j_2<...<j_l$ (et $l$ de m\^eme parit\'e que $k$) : si $r$ est le plus grand indice tel que $i_r>i_{r+1}$, on pousse $e_{i_r}$ successivement vers la droite en utilisant les formules de commutation. On recommence pour les indices plus petits si n\'ec\'essaire.

Ainsi $1$ et les $e_{j_1}\times...\times e_{j_l}$ pour $j_1<j_2<...<j_l$ forment une famille g\'en\'eratrice de l'alg\`ebre de Clifford.
\end{proof}
\section{Calculs explicites en degr\'e $2$}

Calculer $C(A,\s)$ est difficile en g\'en\'eral, m\^eme pour une 
involution. Je n'effectue ici les calculs que pour le plus bas degr\'e 
possible pour montrer~:

\begin{prop}\label{deg 2}
Si $\deg A = 2$, alors $C(A,\s) \simeq F[X]/(X^2-\Nrd(a_\s+1) \disc(\s))$ 
et donc en particulier, si $a_\s-1 \in A^\times$, alors 
$C(A,\s) \simeq F[X]/(X^2-\Nrd(a_\s^2-1))$.
\end{prop}

Bien s\^ur si $\s$ est une involution de type orthogonal, $a_\s=1$ et on retrouve le cas particulier du th\'eor\`eme \ref{structure} : $C(A,\s) \simeq F[X]/(X^2-\disc(\s))$.

La preuve se fait au cas par cas. On \'etudie d'abord le cas o\`u $A$ est 
d\'eploy\'ee en utilisant les formes possibles de $a_\s$. Si $A$ est 
une alg\`ebre de quaternions, on utilise le fait que $\s=(\INT u) \circ \r$, 
pour $u \in A^\times$ et $\r$ l'involution canonique sur $A$ et on 
exprime les calculs en fonction de $u$.

\subsection{Le cas d\'eploy\'e}

On suppose ici que $A=\End V$ et que l'antiautomorphisme $\s=\s_b$,  d'asym\'etrie $a_\s=a$, est adjoint \`a la forme bilin\'eaire $b : V \times V \to F$  non d\'eg\'en\'er\'ee sur $V=F^2$. Dans \cite[Th\'eor\`eme 1]{CT}, on donne des conditions n\'ecessaires et suffisantes pour que $a \in A^\times$ soit une asym\'etrie. En particulier, il 
est n\'ecessaire que $a$ soit conjugu\'e \`a son inverse. Donc 
si  $\l$ est valeur propre de $a$, alors $\l^{-1}$ est aussi valeur 
propre de $a$. Ainsi, pour toute asym\'etrie $a$, il existe 
une base $(i,j)$ de $V$ dans laquelle $a$ a une matrice de 
l'une des formes suivantes~:

$$\begin{pmatrix}
\l & 0 \\
0 & \l^{-1}
\end{pmatrix},\quad
\begin{pmatrix}
1 & 1 \\
0 & 1
\end{pmatrix},\quad
\begin{pmatrix}
-1 & 1 \\
0 & -1
\end{pmatrix},\quad
\begin{pmatrix}
1 & 0 \\
0 & -1
\end{pmatrix},$$
ou une matrice non trigonalisable. Parmi elles, seules les matrices 
$\begin{pmatrix}
\l & 0 \\
0 & \l^{-1}
\end{pmatrix}$ et 
$\begin{pmatrix}
-1 & 1 \\
0 & -1
\end{pmatrix}$ et certaines non trigonalisables v\'erifient les autres 
conditions n\'ecessaires (avec les notations de \cite{CT}, pour 
$\begin{pmatrix}
1 & 1 \\
0 & 1
\end{pmatrix}$, $V_2^1$ est de dimension paire, et pour 
$\begin{pmatrix}
1 & 0 \\
0 & -1
\end{pmatrix}$, $V_1^{-1}$ est de dimension impaire, donc ce 
ne sont pas des asym\'etries). Il reste donc trois cas \`a \'etudier. 

On note ici $I$ l'id\'eal de $T(V)$ engendr\'e par les \'el\'ements du type $a(w)\otimes w-b(w,w)$ pour $w\in V$.

\medskip

\noindent{\underline{\it Premier cas} - } Si $a$ a pour matrice 
$\begin{pmatrix}
\l & 0 \\
0 & \l^{-1}
\end{pmatrix}$ dans la base $(i,j)$, supposons tout d'abord $\l=1$. 
Alors $b$ est bilin\'eaire sym\'etrique et c'est le cas classique. 
Supposons maintenant $\l \not= 1$. Alors, comme, pour tous $x,y \in V$, 
$b(x,y)=b(y,a(x))$, on montre que la matrice de $b$ 
dans $(i,j)$ est, \`a un facteur multiplicatif pr\`es, la matrice $\begin{pmatrix} 0 & 1 \\ \l^{-1} & 0\end{pmatrix}$. 

L'id\'eal $I$ contient ainsi les trois \'el\'ements 
\eqna
a(i) \otimes i - b(i,i)&=&\l i \otimes i,\\
a(j) \otimes j - b(j,j)&=&\l^{-1} j \otimes j \\
\mbox{et}\quad a(i) \otimes j + a(j) \otimes i - b(i,j)-b(j,i)&=&
\l i \otimes j + \l^{-1}j \otimes i-(1+\l^{-1}).
\endeqna
 Donc la partie paire de 
$T(V)/I$ est engendr\'ee par $i \times j$, et comme dans le quotient 
$j \times i=(\l+1) -\l^2 (i \times j)$, on en d\'eduit que 
$(i \times j)^2=(\l+1) (i \times j)$ et donc 
$$C(A,\s)=F[X]/(X(X-(1+\l)))$$ 
qui est isomorphe \`a $F^2$ si $\l \not= -1$ et 
\`a $F[X]/(X^2)$ si $\l=-1$. 

Or $\disc \s = -\det b = \l^{-1}$ et $\Nrd(a+1)=(\l+1)(\l^{-1} +1)$, donc $\Nrd(a+1)\disc(\s)=0$ si $\l=-1$ 
et appartient \`a $(F^\times)^2$ si $\l\neq -1$. D'o\`u le r\'esultat.

\medskip

\noindent{\underline{\it Deuxi\`eme cas} - } 
Si $a$ a pour matrice 
$\begin{pmatrix}
-1 & 1 \\
0 & -1
\end{pmatrix}$ dans la base $(i,j)$,  la matrice 
de $b$ dans $(i,j)$ est (\`a un facteur pr\`es) la matrice 
$\begin{pmatrix} 0 & 1 \\ -1 & -1/2 
\end{pmatrix}$. 

L'id\'eal $I$ contient ainsi $-i \otimes i$, $-j \otimes j + \frac{1}{2}$ et $-i\otimes j - j\otimes i + i\otimes i-1+1$. 
On en d\'eduit que $C(A,\s)$ est engendr\'e par $i \times j$, 
que dans $T(V)/I$, $j \times i = - i \times j$ et donc que $(i \times j)^2 = 0$. 
D'o\`u $C(A,\s)=F[X]/(X^2)$, ce qui est le r\'esultat attendu car 
$\Nrd(a+1)=0$. 

\medskip

\noindent{\underline{\it Troisi\`eme cas} - } 
Si $a$ n'a pas de valeur propre sur $F$. Alors sur la cl\^oture alg\'ebrique 
de $F$, elle a deux valeurs propres distinctes inverses l'une de l'autre 
$\l$ et $\l^{-1}$. Notons $\a=\l+\l ^{-1}$. Alors le polyn\^ome minimal de $a$ est $P=(X-\l)(X-\l^{-1})=X^2-\a X+1=\det(X-a)$. En particulier,  \eqna
\Nrd(a+1)&=\det(a+1)=P(-1)&=2+\a\\
\mbox{et}\quad\disc(\s)&=-\Nrd(1-a)=-P(1)&=-2+\a.
\endeqna 
Finalement $\disc(\s)\Nrd(a+1)=-4+\a^2$ est le discriminant de $P$, et ainsi $$F[X]/(P)=F[X]/(X^2-\disc\s\Nrd(a+1)).$$

De plus, $P$ \'etant irr\'eductible sur $F$, il existe une base $(i,j)$ de $V$ dans laquelle la matrice de $a$ est la matrice compagnon 
$\begin{pmatrix}
0&-1\\1&\a
\end{pmatrix}$ de $P$. On en d\'eduit que la matrice de $b$ est  
$\begin{pmatrix}
-1 & 1 \\
0 & -1
\end{pmatrix}$ \`a un facteur multiplicatif pr\`es.

L'id\'eal $I$ contient alors $j\otimes i-1$, $-i\otimes j+\a j\otimes j-1$ et $j\otimes j-i\otimes i+\a j\otimes i-1-(\a -1)$, donc aussi $j\otimes j-i\otimes i$.
On en d\'eduit que $C(A,\s)$ est engendr\'ee par $j \times j$, 
que dans $T(V)/I$, 
$$j\times j =i\times i,\quad i\times j=\a j\times j-1\quad\mbox{et}\quad j\times i=1$$ et donc que $$(j\times j)^2 = (j\times i)\times(i\times j)=i\times j=\a j\times j-1.$$
Ainsi $j\times j$ a pour polyn\^ome minimal $P$ et   $C(A,\s)=F[X]/(P)$.
Ceci ach\`eve la preuve dans le cas d\'eploy\'e.

\subsection{Pour une alg\`ebre de quaternions} 

Soit $A=(\a,\b)_F$ une alg\`ebre de quaternions sur $F$, dans laquelle on note $(1_A,i,j,k)$ la base usuelle (telle que $i^2=\a, j^2=\b$ et $ij=-ji=k$). Soit $\r$ l'involution canonique (de type symplectique, $\r~:~A\rightarrow A\,;\,q=x+yi+zj+tk\mapsto x-yi-zj-tk$).

Rappellons que d'apr\`es le th\'eor\`eme de Skolem-Noether, si $\s$ est un antiautomorphisme de $A$, il existe $u\in A^\times$ tel que $\s=\INT u\circ\r$. D'apr\`es \cite[proposition 7]{CT}, et puisque $a_\r=-1$, cela donne $a_\s=u\r(u)^{-1}a_r=-u\r(u)^{-1}$.

Notons enfin que si $q\in A$ alors $\Trd q$ est la coordonn\'ee de $q+\r(q)$ sur $1_A$.

Nous distinguons trois cas.

\noindent{\underline{\it Premier cas} - } si $\s$ est une involution orthogonale, il est connu que $C(A,\s)=F(u)$. Or $1=a_\s$ donc $u$ est un quaternion pur, et quitte \`a changer de base, on peut supposer que $u=i$ et donc que $\s(q)=x-yi+zj+tk$. Ainsi, $C(A,\s)=F(i)=F[X]/(X^2-\a)$, et comme $i\in\Skew(A,\g_\s)\bigcap A^\times$, $\Nrd(a_\s+1)\disc\s=-\Nrd 2\Nrd i=4\a$, on obtient le r\'esultat souhait\'e.

\medskip

\noindent{\underline{\it Deuxi\`eme cas} - } Si $\s=\r$ (la seule involution symplectique).

Comme $\g_\r(q)=-\r(q)$, $q\in A$ est sym\'etrique pour $\g_\r$ si et seulement si sa trace est nulle et $J_1$ est engendr\'e par les quaternions purs, \cad par $i,j,k$.

De plus 
\eqna
s\in\Sym(\Aun\otimes\Aun,\g_{\rhot ,2})&\Leftrightarrow & \forall q\in A \quad \Sand s(-\r(q))=\Sand s(q)\\
&\Leftrightarrow &\forall q\in A \quad \Sand s(q+\r(q))=0\\
&\Leftrightarrow & \Sand s(1)=\mu(s)=0.
\endeqna
De plus $\mu=-\mu_\s$, donc $J_2$ est engendr\'e par les $s+\frac{1}{2}\mu(s)$ pour $\mu(s)=0$, et n'est autre que l'id\'eal engendr\'e par $\ker \mu$.

Ainsi l'alg\`ebre $C(A,\r)= T(\Aun)/(J_1 + J_2)$ est engendr\'ee par $1_A$, et comme $$\mu(\a\b1_A\otimes 1_A+k\otimes k)=\a\b1_A.1_A+k.k=0,$$ on obtient $1_A^2=0$ dans $C(A,\r)$. Ceci signifie que cette alg\`ebre est exactement $$C(A,\r)=F[X]/(X^2)=F[X]/(X^2-\Nrd(a+1)\disc\r).$$

\medskip

\noindent{\underline{\it Troisi\`eme cas} - } Si maintenant $\s$ n'est pas involutive, \cad si $a_\s\neq\pm 1$, alors $u$ n'est ni un scalaire, ni un quaternion pur. 
Donc $u$ est de la forme $\l(1+i_0)$ o\`u $i_0$ est un quaternion pur. Ainsi $\INT u=\INT(1+i_0)$ et on peut donc se ramener \`a une base $(1_A,i,j,k)$ de $A$ telle que $u=1+i_0$.

Remarquons que l'on peut \'ecrire tout \'el\'ement de $A$ de mani\`ere unique sous la forme $q=v+wj$ avec $v,w\in F(i)$, et que ce sous-corps commutatif $F(i)$ de $A$ est stable par $\r$ et $\s$ et contient $u$ et donc aussi $a_\s$. De plus si $v\in F(i)$ alors $vj=j\r(v)$. 

Nous pouvons maintenant \'ecrire :
\eqna\displaystyle
a_\s&=&-u\r(u)^{-1}=-\frac{1+i}{1-i},\\ \g_\s(q)&=&\s(q)a_\s=-(1+i)\r(q)(1+i)^{-1}\frac{1+i}{1-i}=-(1+i)\r(q)(1-i)^{-1} \\  \rm{et}\quad\g_{\sigt}(q)&=&a_\s\s(q)a_\s^2=-\frac{(1+i)^2}{1-i}\r(q)\frac{1+i}{(1-i)^2},
\endeqna
 \cad si $q=x+yi+zj+tk=(x+yi)+(z+ti)j$, alors

$$\displaystyle{\g_\s(q)=(-x+yi)\frac{1+i}{1-i}+(z+ti)j\quad\mbox{et}
\quad\g_\sigt(q)=(-x+yi)\Bigl(\frac{1+i}{1-i}\Bigr)^3+(z+ti)j}.$$

Ainsi
\eqna
&\g_{\s}(q)=q\\
\Leftrightarrow& (-x+yi)(1+i)=(x+yi)(1-i)\\
\Leftrightarrow& x=\a y.
\endeqna
  
On en d\'eduit que $\Sym(a,\g_\s)=<\a 1_A+i,j,k>$ et donc 

$$J_1=<\a 1_A+i-\frac{1}{2}\Trd(\a 1_A+i),j-\frac{1}{2}\Trd j,k-\frac{1}{2}\Trd k>=<i+\a(1_A-1_F),j,k>.$$

De plus 
\eqna
&\g_{\sigt}(q)=-q&\\
\Leftrightarrow& z=t=0\quad\rm{et}&\displaystyle{(-x+yi)\Bigl(\frac{1+i}{1-i}\Bigr)^3}=-x-yi\\
\Leftrightarrow& z=t=0\quad\rm{et}& y(1+3\a)=x(\a+3).
\endeqna
Alors comme $1+3\a$ et $\a+3$ ne sont pas conjointement nuls, $\Skew(A,\g_{\sigt})$ est un sous espace vectoriel de $A$ de dimension $1$. Or il contient $(1+i)^3$, qui en est donc un g\'en\'erateur.

Mais $s\in\Sym(A,\g_{\sigt,2})\Leftrightarrow \forall x\in A\quad(\Sand s)(\g_{\sigt}(x)-x)=0$, et comme $\g_{\sigt}$ est une une application lin\'eaire involutive, cela \'equivaut encore \`a $\forall x\in\Skew(A,\g_{\sigt})\quad (\Sand s)(x)=0$ et donc \`a $$(\Sand s)(1+i)^3=0.$$

Finalement, dans $C(A,\s)$, $j=k=0$ et $i=\a(1_F-1_A)$, ce qui signifie que $C(A,\s)$ est engendr\'ee par $1_A$ (ou par $i$). Mais comme $1_A$ et $i$ commutent dans $A$ avec $1+i$ et donc avec $(1+i)^3$ et $a$, on obtient 
\eqna
(\Sand(\a 1_A\otimes 1_A-i\otimes i))(1+i)^3&=&\a(1+i)^3-i^2(1+i)^3=0\quad\mbox{et} \\
\mu_\s(\a 1_A\otimes 1_A-i\otimes i)&=&-\a a+i^2a=0.
\endeqna
On en d\'eduit que l'\'el\'ement $\a 1_A\otimes 1_A-i\otimes i-\frac{1}{2}\mu(\a 1_A\otimes 1_A-i\otimes i)$ est dans $J_2$ et donc que dans $C(A,\s)$, $$\a 1_A\times 1_A=i\times i=(\a(1_F-1_A))^2.$$

Ceci prouve que $1_A$ a pour polyn\^ome minimal $(\a-1)X^2-2\a X+\a$, dont le discriminant modulo les carr\'es est $\a$, d'o\`u $$C(A,\s)=F[X]/(X^2-\a)$$ (et est isomorphe, comme dans le cas d'une involution orthogonale au sous-corps $F(u)$ de $A$).

De plus, ici $a_\s-1\neq 0$ donc $\disc\s=-\Nrd(a_\s-1)$ et ainsi $$\displaystyle\disc\s\Nrd(a_\s+1)=-\Nrd(a_\s^2-1)=\Nrd\Bigl(\frac{4i}{(1-i)^2}\Bigr)=\a$$ dans $F^\times/F^{\times 2}$, ce qui ach\`eve la d\'emonstration.

{\sc Remerciements} : Je remercie chaleureusement Anne Queguiner-Mathieu pour ses remarques et conseils judicieux sur ce travail et pour l'int\'er\^et qu'elle y a port\'e.

\end{document}